\documentclass[reqno]{amsart}
\usepackage{graphicx,subfigure}

\numberwithin{equation}{section}
\usepackage[utf8]{inputenc}
\usepackage[english]{babel}

\usepackage{amsmath,amsthm,amsfonts,latexsym,amssymb}
\usepackage[colorlinks]{hyperref}
\hypersetup{linkcolor=blue,citecolor=blue,filecolor=black,urlcolor=blue}
\usepackage{comment}

\usepackage{color,amsthm,amsfonts}

{ \theoremstyle{plain}
\newtheorem{theorem}{Theorem}[section]
\newtheorem{proposition}[theorem]{Proposition}
\newtheorem{lemma}[theorem]{Lemma}

  \theoremstyle{remark}
\newtheorem{remark}[theorem]{Remark}
  \theoremstyle{definition}

}

\def\RR{\mathbb{R}}
\def\NN{\mathbb{N}}

\begin{document}

\title[A $p$-Laplacian supercritical Neumann problem]{A \MakeLowercase{$p$}-Laplacian Neumann problem with a possibly supercritical nonlinearity}\author[F. Colasuonno]{Francesca Colasuonno}

\maketitle

\begin{abstract}
We look for nonconstant, positive, radially nondecreasing solutions of the quasilinear equation $-\Delta_p u+u^{p-1}=f(u)$ with $p>2$, in the unit ball $B$ of $\mathbb R^N$, subject to homogeneous Neumann boundary conditions. The assumptions on the nonlinearity $f$ are very mild and allow it to be possibly supercritical in the sense of Sobolev embeddings. The main tools used are the truncation method and a mountain pass-type argument. In the pure power case, i.e., $f(u)=u^{q-1}$, we detect the  limit profile of the solutions of the problems as $q\to\infty$. 
\end{abstract}

\section{Introduction and main results}
In \cite{CN}, we study the existence of nonconstant, radially nondecreasing solutions of the following quasilinear problem
\begin{equation}\label{P}
\begin{cases}
-\Delta_p u+u^{p-1}=f(u)\quad&\mbox{in }B,\\
u>0&\mbox{in }B,\\
\partial_\nu u=0&\mbox{on }\partial B,
\end{cases}
\end{equation}
where $B$ is the unit ball of $\RR^N$, $N\ge1$, $\nu$ is the outer unit normal of $\partial	B$, and $\Delta_p u:=\mathrm{div}(|\nabla u|^{p-2}\nabla u)$ is the $p$-Laplacian operator, with $p>2$. 
We require very mild assumptions on the nonlinearity $f$ on the right-hand side, namely $f\in C^1([0,\infty))$ and satisfies the following hypotheses

($f_1$) $\lim_{s\to 0^+}\frac{f(s)}{s^{p-1}}\in [0,1)$;

($f_2$) $\liminf_{s\to+\infty}\frac{f(s)}{s^{p-1}}>1$;

($f_3$) $\exists$ a constant $u_0>0$ such that $f(u_0)=u_0^{p-1}$ and $f'(u_0)>(p-1)u_0^{p-2}$.	
\medskip 

Our main results in \cite{CN} read as follows.

\begin{theorem}\label{main}
If $f$ satisfies \rm{($f_1$)-($f_3$)}, there exists a nonconstant, radially nondecreasing solution of \eqref{P}. If furthermore there exist $n$ different positive constants $u_0^{(1)}\neq\dots\neq u_0^{(n)}$ for which \rm{($f_3$)} holds, then \eqref{P} admits at least $n$ distinct nonconstant, radially nondecreasing solutions.
\end{theorem}

\begin{theorem}\label{limit_profile}
Let $f(u)=u^{q-1}$, with $q>p$. Denote by $u_q$ the solution found in Theorem \ref{main}, corresponding to such $f$. Then, as $q\to\infty$, 
$$u_q\to G\mbox{ in }W^{1,p}(B)\cap C^{0,\mu}(\overline{B})\quad\mbox{for any }\mu\in(0,1),$$
where $G$ is the unique solution of the Dirichlet problem
$$
\begin{cases}
-\Delta_p G+G^{p-1}=0\quad&\mbox{in }B,\\
G=1&\mbox{on }\partial B.
\end{cases} 
$$
\end{theorem}

{\sc Remarks.}

\noindent $\bullet$ We observe that $f$ is allowed to be supercritical in the sense of Sobolev embeddings, which will be the most interesting case.\smallskip

\noindent $\bullet$ The model $f$ is the pure power function $f(u)=u^{q-1}$, with $q>p$. In this case, problem \eqref{P} admits the constant solution $u\equiv 1$ for every $q>p$, including the supercritical  case $q>p^*$, where $p^*:=Np/(N-p)$ if $p<N$ and $p^*:=+\infty$ otherwise.
Therefore, the natural question that arises is whether \eqref{P} admits any {\it nonconstant} solutions.
It is worth stressing a remarkable difference between problem \eqref{P} and the analogous problem under homogeneous Dirichlet boundary conditions. Indeed, it is well-known that, as a consequence of the Poho\v{z}aev identity (cf. \cite[Section 2]{PS}), the Dirichlet problem does not admit any nonzero solutions when $q\ge p^*$.

\noindent $\bullet$ We remark that condition ($f_3$) is absolutely natural under ($f_1$) and ($f_2$). Indeed, by the regularity of $f$ and by ($f_1$)-($f_2$), there must exist an intersection point $u_0$ between $f$ and the power $s^{p-1}$ such that $f'(u_0)\ge (s^{p-1})'(u_0)=(p-1)u_0^{p-2}$. Hence, ($f_3$) is only meant to exclude the possibility of a degenerate situation in which $f$ is tangent to $s^{p-1}$ at $u_0$. \smallskip

\noindent $\bullet$ We can always think $f$ to satisfy also 
\smallskip

($f_0$)\quad $f\ge0$ and $f'\ge0$.
\smallskip

\noindent Indeed, if this is not the case, we can replace $f$ by $g(s):=f(s)+(m-1)s^{p-1}$ for a suitable $m>1$ such that $g\ge0$ and $g'\ge 0$, and study the equivalent problem 
$$
\begin{cases}
-\Delta_p u+m u^{p-1}=g(u)\quad&\mbox{in }B,\\
u>0&\mbox{in }B,\\
\partial_\nu u=0&\mbox{on }B.
\end{cases}
$$
Therefore, without loss of generality, {\it from now on in the paper we assume $f$ to satisfy \rm{($f_0$)} as well.} 

\begin{figure}
\centering
\includegraphics[scale=0.36]{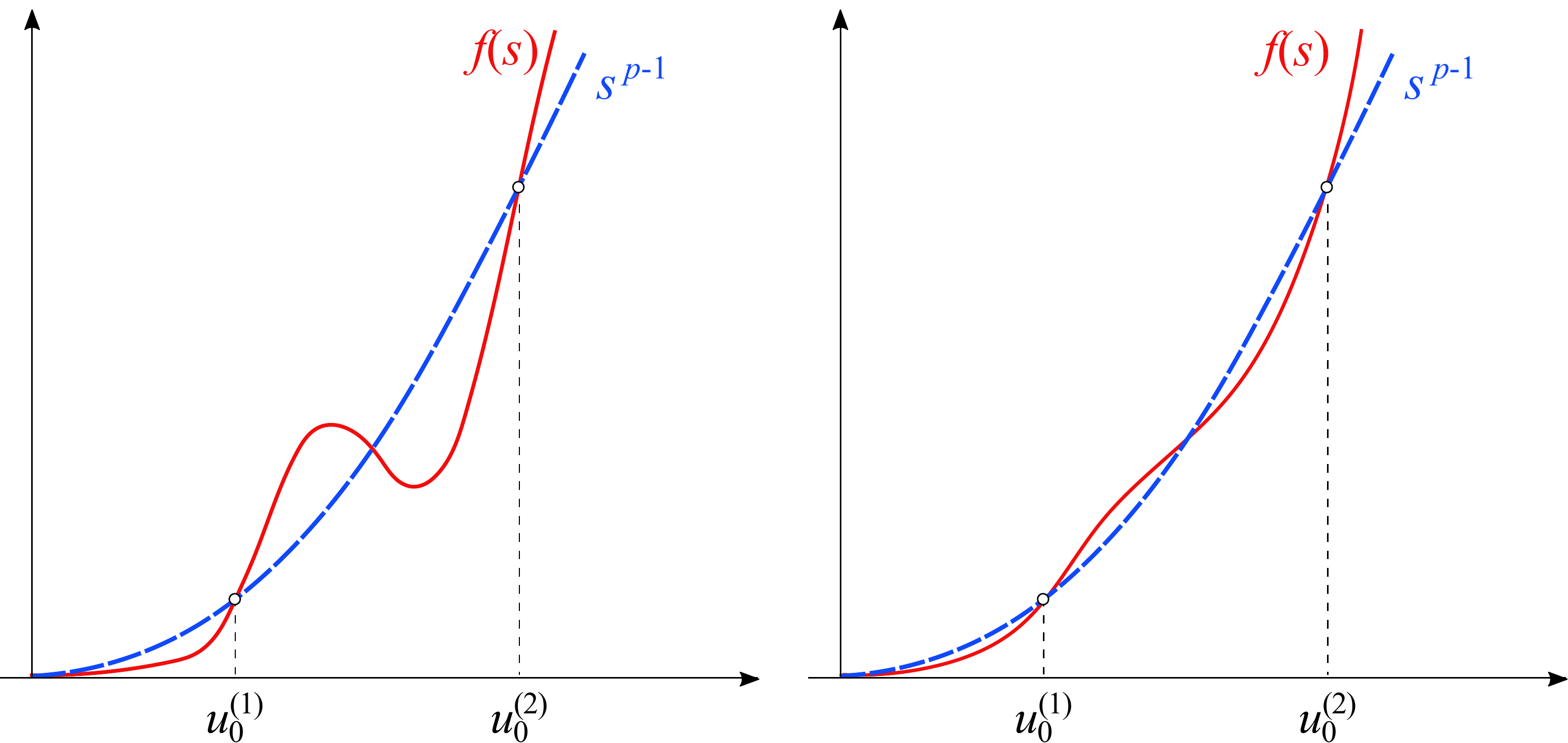}
\caption{{\small {\it Left}: Graph of a sample nonlinearity $f$ satisfying ($f_1$)-($f_3$). {\it Right}: Graph of a sample nonlinearity $f$ satisfying ($f_0$)-($f_3$).}}
\end{figure}

Since $f$ is possibly supercritical, the energy functional $I$ associated to the problem is not well defined in the whole of $W^{1,p}(B)$, and so a priori we cannot use variational techniques to solve the problem. This issue is overcome for the first time in \cite{ST} for the semilinear case ($p=2$) and then in \cite{S} for any $1<p<\infty$, by working in the closed and convex cone
$$\mathcal C :=\left\{u\in W^{1,p}_{\mathrm{rad}}(B)\,:\, u\ge 0\;\mbox{ and }\;u(r)\le u(s) \mbox{ for }r\le s\right\},$$
where we have denoted by $W^{1,p}_{\mathrm{rad}}(B)$ the space of $W^{1,p}(B)$-functions which are radially symmetric and with abuse of notation we have written $u(x)=u(r)$ for $|x|=r$.  
Indeed, this cone has the property that all its functions are bounded, i.e.,
\begin{equation}\label{boundedfunct}
\|u\|_{L^\infty(B)}\le C(N)\|u\|_{W^{1,p}(B)}\quad\mbox{for some $C(N)>0$ independent of $u\in\mathcal C$},
\end{equation}
see e.g. \cite[Lemma 2.2]{CN}. Due to \eqref{boundedfunct}, it makes sense to define an energy functional $I$ in $\mathcal C$, associated to the equation. On the other hand, the main disadvantage for working in this cone is the fact that it has empty interior in the $W^{1,p}$-topology. As a consequence, in general, critical points of $I:\mathcal C\to\mathbb R$ are not solutions of \eqref{P}. In \cite{ST,S}, the authors require additional assumptions on $f$ to prove that the critical point of $I$, found via variational techniques, is indeed a weak solution of the problem. While in \cite{BNW}, in order to weaken the hypotheses on $f$, a different strategy based on the truncation method is proposed. 

The techniques that we use in \cite{CN} to prove Theorem \ref{main} are essentially in the spirit of \cite{BNW}. The scheme of the proof can be split into five steps. 

\underline{{\sc Step 1.}} We first obtain, in \cite[Lemma 2.5]{CN}, the following a priori estimate 
$$\|u\|_{L^\infty(B)}\le K_\infty\quad\mbox{ for all }u\in\mathcal C \mbox{\; that solves \;}\eqref{P},$$ 
for some $K_\infty>0$ independent of $u$. Clearly, $K_\infty\ge u_0$, being $u\equiv u_0$ a solution of \eqref{P} belonging to $\mathcal C$. \smallskip

\underline{{\sc Step 2.}}  This allows us to truncate the nonlinearity $f$, in order to deal with a subcritical nonlinearity $\tilde f$, and so in \cite[Lemma 3.1]{CN}, we prove that 
\begin{center}
{\it For all  $\ell\in(p,p^*)$ there exists $\tilde f\in C^1([0,\infty))$ satisfying \rm{($f_0$)-($f_3$)}, 
\vspace{-.2cm} $$\lim_{s\to\infty}\frac{\tilde f(s)}{s^{\ell-1}}=1,\quad\mbox{ and }\quad\tilde f=f\mbox{\; in \;}[0,K_\infty].$$}
\end{center}

We introduce the following auxiliary problem 
\begin{equation}\label{tildeP}
\begin{cases}
-\Delta_p u+u^{p-1}=\tilde f(u)\quad&\mbox{in }B,\\
u>0&\mbox{in }B,\\
\partial_\nu u=0&\mbox{on }\partial B.
\end{cases}
\end{equation}
As a consequence of the previous two steps, it is immediate to see that 
\begin{center}
{\it In the cone $\mathcal C$, the two problems \eqref{P} and \eqref{tildeP} are equivalent.}
\end{center}

\underline{{\sc Step 3.}} Thanks to the subcriticality of $\tilde f$, we can define the energy functional associated to \eqref{tildeP} in the whole of $W^{1,p}(B)$ as follows
$$\tilde I(u):=\frac1p\int_B(|\nabla u|^p+|u|^p)dx-\int_B\tilde F(u)dx,\quad\mbox{where }\tilde F(u):=\int_0^u \tilde f(s)ds$$
for all $u\in W^{1,p}(B)$. All critical points of $\tilde I$ are weak solutions of \eqref{tildeP}.

\begin{remark} 
Since $p>2$, $\tilde I$ is of class $C^2$, while if $1<p<2$, the functional $\tilde I$ is only of class $C^1$. This lack of regularity prevents either the use of second order Taylor expansions as done in \cite{BNW,CN} (see also Section \ref{sec3} below) or the use of a generalized Morse Lemma when looking for {\it nonconstant} solutions. Moreover, when $1<p<2$, Simon's inequalites relating $\tilde I'$ and the pseudo-differential gradient are weaker than the ones found for the case $p>2$, this makes harder the construction of a descending flow and consequently the proof of a deformation-type lemma. 
\end{remark}

\underline{{\sc Step 4.}} We find a critical point $u$ of $\tilde I$ belonging to $\mathcal C$ via a mountain pass-type argument. We localize the solution in such a way that, if we have $n$ different positive constants $u_0^{(i)}$ verifying ($f_3$), we get ``for free'' also the multiplicity result stated in Theorem \ref{main}.\smallskip

\underline{{\sc Step 5.}} We prove that the solution found in Step 4. is nonconstant, by using a second order Taylor expansion of $\tilde I$.
\medskip

In the next two sections we give some details about Steps 4. and 5., respectively. 
While in the last section we sketch the proof of Theorem \ref{limit_profile}. 

\section{Step 4: A nonconstant solution of \eqref{P} belonging to $\mathcal C$}\label{sec2}

Due to the subcriticality of $\tilde f$, it is standard to prove the following compactness result (cf. \cite[Lemma 3.4]{CN}):
\begin{center}
{\it The functional $\tilde I$ satisfies the Palais-Smale condition.}   
\end{center}
\medskip

\noindent\textbf{The restricted cone $\mathcal C_*$.}\bigskip

\begin{minipage}[c]{.4\textwidth}
\hspace{-.2cm}\includegraphics[scale=0.38]{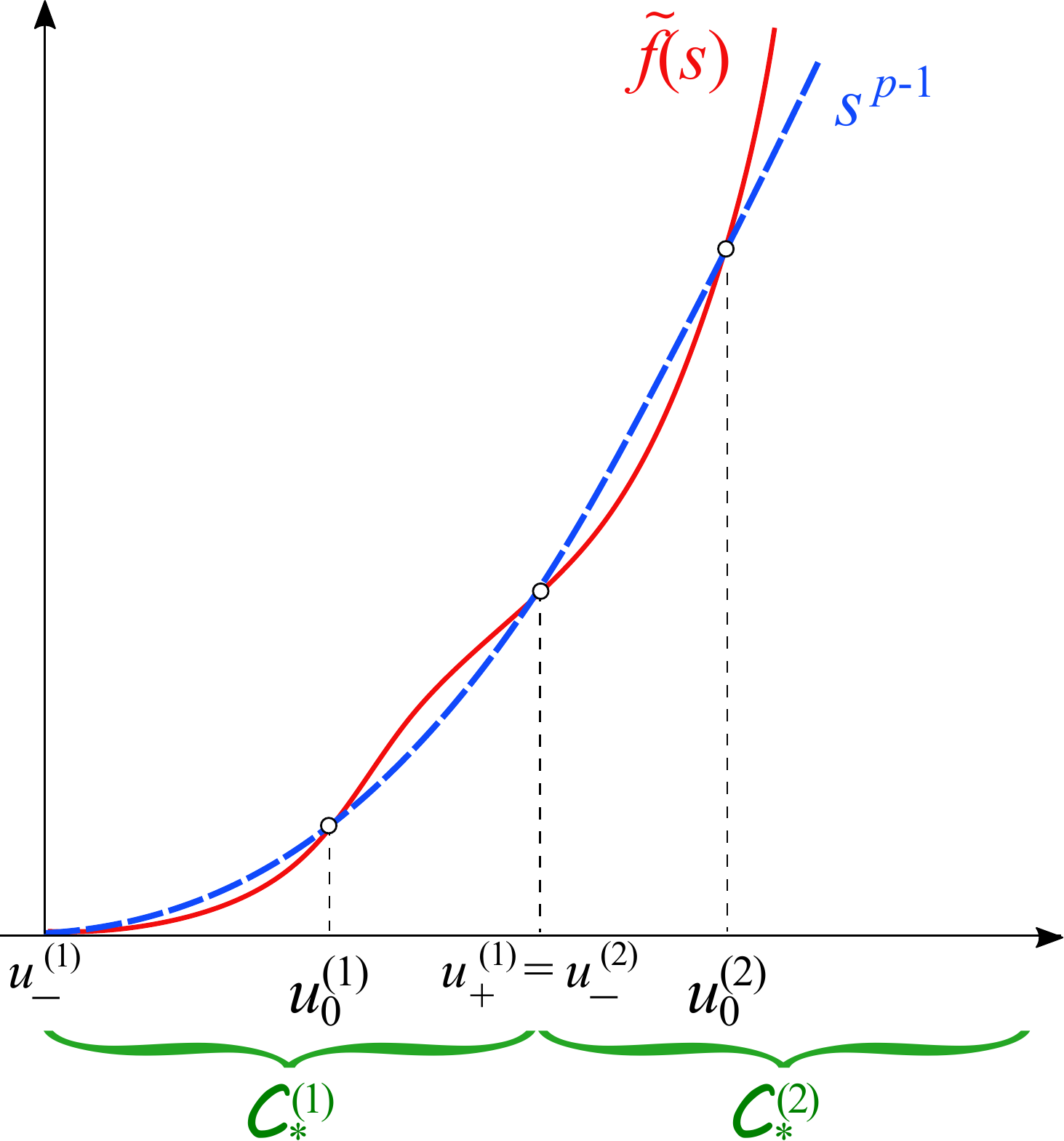}
\end{minipage}\;\qquad
\begin{minipage}[c]{.5\textwidth}
Let $n\in \NN$ be the number of positive constants $u_0^{(i)}$ satisfying ($f_3$). For every $i=1,\dots,n$, we set 
$$
\begin{aligned}
u_-^{(i)}&:=\sup\left\{s\in[0,u_0^{(i)})\,:\, \tilde f(s)=s^{p-1}\right\},\\
u_+^{(i)}&:=\inf\left\{s\in(u_0^{(i)},\infty)\,:\,\tilde f(s)=s^{p-1}\right\}.
\end{aligned}
$$
For every $i$, we further introduce the following subset of $\mathcal C$
$$\mathcal C_*^{(i)}:=\left\{u\in\mathcal C\,:\, u_-^{(i)}\le u\le u_+^{(i)} \right\}$$
which turns out to be itself a closed convex cone of $W^{1,p}(B)$.
\end{minipage}
\medskip

{\sc Remarks. }

\noindent $\bullet$ Thanks to ($f_3$), each $u_0^{(i)}$ is an isolated zero of $\tilde f(s)-s^{p-1}$, hence $u_-^{(i)}\neq u_0^{(i)}\neq u_+^{(i)}$ for every $i=1,\dots,n$.

\noindent $\bullet$ We observe that $u_+^{(n)}$ can be possibly $+\infty$. For instance, for the pure power function $f(u)=u^{q-1}$ with $q>p$, it results $n=1$, $u_-=0$, $u_0=1$, $u_+=+\infty$, and $\mathcal C=\mathcal C_*$.

\noindent $\bullet$ All and only the zeros of $\tilde f(s)-s^{p-1}$ are constant solutions of \eqref{tildeP}, and so of \eqref{P}. Hence, the only constant solutions of \eqref{P} belonging to $\mathcal C_*^{(i)}$ are $u_-^{(i)},\,u_0^{(i)}$, and $u_+^{(i)}$.  

\noindent $\bullet$ If we prove the existence of a nonconstant solution $u$ belonging to $\mathcal C_*^{(i)}$, we know at once that $u_-^{(i)}\le u\le u_+^{(i)}$ and that $u_-^{(i)}\not\equiv u\not\equiv u_+^{(i)}$. This implies that nonconstant solutions of \eqref{P} belonging to different $\mathcal C_*^{(i)}\,$'s are different.
\smallskip

As a consequence of the last two remarks, we can see that the advantage of working in $\mathcal C_*^{(i)}$ instead of $\mathcal C$ is twofold. Firstly, it helps avoiding constant solutions: it is enough to prove that the solution found is none of the three constant solutions in $\mathcal C_*^{(i)}$. Secondly, the restricted cone $\mathcal C_*^{(i)}$ allows us to localize our solution, so that the multiplicity part of Theorem \ref{main} follows immediately by the existence part.  
   
{\it Hereafter, we assume for simplicity $n=1$ and we omit all the superscripts $(i)$.} Clearly, if $n>1$, it is possible to repeat the same arguments in each cone $\mathcal C_*^{(i)}$.     
   
\subsection*{A deformation lemma}
This is the most technical part of the proof. Since the space $W^{1,p}(B)$ in which the energy functional $\tilde I$ is defined is bigger than the set $\mathcal C_*$ in which we want to find a minimax solution, we need a slightly different version of the deformation lemma. 

\begin{lemma}[Lemma 3.9 of \cite{CN}]\label{deformation} Let $c\in\mathbb R$ be such that $\tilde I'(u)\neq 0$ for all $u\in \mathcal C_*$, with $\tilde I(u)=c$. Then, there exist a positive constant $\bar\varepsilon$ and a function $\eta:\mathcal C_*\to\mathcal C_*$ satisfying the following properties: 
\begin{itemize}
\item[(i)] $\eta$ is continuous with respect to the topology of $W^{1,p}(B)$;
\item[(ii)] $I(\eta(u))\le I(u)$ for all $u\in\mathcal C_*$;
\item[(iii)] $I(\eta(u))\le c-\bar\varepsilon$ for all $u\in\mathcal C_*$ such that $|I(u)-c|<\bar\varepsilon$;
\item[(iv)] $\eta(u)=u$ for all $u\in\mathcal C_*$ such that $|I(u)-c|>2\bar\varepsilon$.
\end{itemize}
\end{lemma}
{\sc Remarks.}

\noindent $\bullet$ We stress here that we build a deformation $\eta$ not only for regular values $c$ of $\tilde I$ (i.e., such that  $\tilde I'(u)\neq 0$ \textit{for all }$u\in W^{1,p}(B)$ with $\tilde I(u)=c$), but also for all $c\in \RR$ for which $\tilde I'(u)\neq 0$ \textit{for all }$u\in\mathcal C_*$ with $\tilde I(u)=c$.

\noindent $\bullet$ In this version of the deformation lemma, we need to prove that the $\eta$ preserves the cone $\mathcal C_*$. This is the most delicate point of the proof. It requires the existence of a pseudo-gradient vector field $K$ of $\tilde I$ which is not only locally Lipschitz continuous, but which satisfies also the following property
\begin{equation}\label{Kpreserves}
K\left(\mathcal C_*\setminus\{\mbox{critical points of $\tilde{I}$}\}\right)\subset\mathcal C_*.
\end{equation}
Indeed, for every $u\in\mathcal C_*$, the deformation $\eta(u)$ is built as the unique solution $\mu(t,u)$ of the Cauchy problem 
\begin{equation}\label{pdc}
\begin{gathered}
\begin{cases}
\frac{d}{dt}\mu(t,u(x))=-\Phi(\mu(t,u(x)))\quad&(t,x)\in(0,\infty)\times B,\\
\partial_\nu\mu(t,u(x))=0&(t,x)\in(0,\infty)\times\partial B,\\
\mu(0,u(x))=u(x) & x\in B,
\end{cases}\\
\mbox{where }
\Phi(u):=\begin{cases}\chi_1(I(u))\chi_2(u)\frac{u-K(u)}{\|u-K(u)\|}\,&\mbox{if }|I(u)-c|\le 2\bar\varepsilon,\\
0&\mbox{otherwise,}\end{cases}\quad\chi_1,\,\chi_2\mbox{ cutoff}
\end{gathered}
\end{equation}
for $t$ (fixed) sufficiently large (i.e., $\eta(u):=\mu(\bar t,u)$). The existence of such operator $K$ and of its properties are proved in \cite[Proposition 3.2 and Lemmas 3.5-3.8]{CN} (see also \cite{BLW} for the case of an open cone) and passes through the study of an auxiliary operator $\tilde T$ related to the inverse of $-\Delta_p(\cdot)+|\cdot|^{p-2}(\cdot)$. In particular, property \eqref{Kpreserves} is a consequence of the fact that $\tilde T(\mathcal C_*)\subseteq \mathcal C_*$, that is proved --by hands-- in \cite[Lemma~3.5]{CN}.
Finally, thanks to \eqref{Kpreserves}, the convexity, and the closedness of $\mathcal C_*$, we are able to prove that $\eta(\mathcal C_*)\subseteq\mathcal C_*$.

\noindent $\bullet$ Condition $(iv)$ is an immediate consequence of the fact that $\mu$ solves the Cauchy problem \eqref{pdc}. While, $(ii)$ and $(iii)$ rely essentially on Simon-type inequalities, that is to say relations between $\tilde I'$ and $K$, see \cite[Proposition 3.2 and Lemmas 3.6-3.8]{CN}. 
 
\subsection*{A mountain pass-type geometry}
\begin{lemma}[Lemma 3.10 and formula (3.32) of \cite{CN}]\label{mountainpassgeo}
Let $\tau$ be a constant such that $0<\tau <\min \{u_0-u_-,u_+-u_0\}$. 
Then there exists $\alpha>0$ such that 
\begin{itemize}
\item[(i)] $\tilde I(u)\ge \tilde I(u_-)+\alpha$ for every $u \in
\mathcal C_*$ with $\|u-u_-\|_{L^\infty(B)}=\tau$;
\item[(ii)] if $u_+< \infty$, then $\tilde I(u)\ge \tilde I(u_+)+\alpha$ for every $u \in
\mathcal C_*$ with $\|u-u_+\|_{L^\infty(B)}= \tau$.
\end{itemize}
Furthermore, 
\begin{itemize}
\item[(iii)] $\tilde I(t\cdot 1)\to -\infty$ as $t\to+\infty$.
\end{itemize}
\end{lemma}
{\sc Remarks.}

\noindent $\bullet$ If $u_+=+\infty$, then $(i)$ and $(iii)$ are pretty much the classical conditions required for the mountain pass geometry centered at $u_-$. 

\noindent $\bullet$ If $u_+<+\infty$, then the roles played by $u_-$ and $u_+$ are interchangeable, hence we prove that the points on the sphere $\partial B_{\tau}(u_-):=\{u\in\mathcal C_*\,:\;\|u-u_-\|_{L^\infty(B)}=\tau\}$ and those on $\partial B_{\tau}(u_+):=\{u\in\mathcal C_*\,:\;\|u-u_+\|_{L^\infty(B)}=\tau\}$ satisfy the same condition with respect to $u_-$ and to $u_+$, respectively. In this case, since $u_0-u_->\tau$ and $u_+-u_0>\tau$, 
then the two closed balls $\overline B_{\tau}(u_-)$ and $\overline B_{\tau}(u_+)$ are disjoint. Therefore, suppose --to fix ideas-- that $\tilde I(u_-)\le \tilde I(u_+)$. By $(ii)$, for all $u\in\partial B_\tau(u_+)$ it results $\tilde I(u)\ge \tilde I(u_+)+\alpha$ and there exists $u_-$, for which 
$$\|u_--u_+\|_{L^\infty(B)}>\tau\quad\mbox{and}\quad\tilde I(u_-)\le \tilde I(u_+).$$

\noindent$\bullet$ We remark that in $(i)$ and $(ii)$ it is possible to use the $L^\infty$-norm instead of the $W^{1,p}$-norm, because $\mathcal C_*$-functions are bounded by \eqref{boundedfunct}. In particular, the use of the $L^\infty$-norm allows us to simplify the constants.

\subsection*{Existence of a solution of \eqref{P} in $\mathcal C_*$}
Let $\tau$ and $\alpha$ be the constants introduced in the previous subsection, 
$$
\begin{aligned} U_- &:= \left\{u \in \mathcal C_* \::\: \tilde I(u)<\tilde I(u_-)+\frac{\alpha}{2},\:
\|u-u_- \|_{L^\infty(B)} < \tau\right\},\\
U_+&:=\begin{cases}
\displaystyle{\left\{u \in \mathcal C_* \::\: \tilde I(u)<\tilde I(u_+)+\frac{\alpha}{2},\:
\|u-u_+ \|_{L^\infty(B)} < \tau\right\}},&\mbox{ if }u_+<\infty,\\
&\\
\left\{u \in \mathcal C_* \, :\, \tilde I(u)< \tilde I(u_-),\, \|u-u_-\|_{L^\infty(B)}>\tau\right\},&\mbox{ if }u_+=\infty
\end{cases}
\end{aligned}
$$
the sets from/to which the admissible paths used to define the minimax level start/arrive, 
$$
\Gamma:=\left\{ \gamma\in C([0,1];\mathcal C_*)\ :\  \gamma(0) \in U_-,\:
  \gamma(1) \in U_+\right\}
$$
the set of admissible paths, and 
\begin{equation}\label{minmax}
c:=\inf_{\gamma\in\Gamma}\max_{t\in[0,1]} \tilde I(\gamma(t))
\end{equation}
the minimax level.

By combining together the compactness condition, the mountain pass-type geometry of $\tilde I$, and the deformation lemma presented above, we are able to prove the following result. 

\begin{proposition}[Proposition 3.11 of \cite{CN}]\label{mountainpass} 
The value $c$ defined in \eqref{minmax} is finite and there exists a
critical point $u\in\mathcal C_*\setminus\{u_-,u_+\}$ of $\tilde I$ such that $\tilde I(u)=c$ and $u>0$. In particular, $u$ is a weak solution of \eqref{P}.
\end{proposition}

{\sc Remarks.}

\noindent$\bullet$ We observe that, since every admissible path $\gamma\in\Gamma$ starts from $B_\tau(u_-)$ and arrives in $B_\tau(u_+)$,  due to its continuity, it must cross the sphere $\partial B_\tau(u_-)$ (and also $\partial B_\tau(u_+)$ if $u_+<+\infty$). Then, by  Lemma \ref{mountainpassgeo}-$(i)$ (and also by $(ii)$ if $u_+<+\infty$), 
$$\tilde I(u_-)<c<+\infty\quad\mbox{(and also $\tilde I(u_+)<c$ if $u_+<+\infty$)}.$$
This immediately excludes the possibility that the solution $u$ is the constant $u_-$ (or the constant $u_+$ when this latter is finite).

\noindent$\bullet$ By the maximum principle \cite[Theorem 5]{V}, $u$ is positive. 

\section{Step 5: The solution found is nonconstant.}\label{sec3}
In this section we conclude the proof of Theorem \ref{main}. As already observed in Section \ref{sec2}, the multiplicity part of the theorem follows easily when one works in the restricted cone $\mathcal C_*$. Concerning the nonconstancy of the solution, we already know by Proposition~\ref{mountainpass} that the solution $u\in\mathcal C_*$, at level $c$, is neither the constant $u_-$ nor the constant $u_+$. It remains to show that $u\not\equiv u_0$. In particular, we prove that $c=\tilde I(u)<\tilde I(u_0)$. 

By the very definition of $c$, it is enough to find an admissible path $\bar\gamma$ such that 
\begin{equation}\label{toprove}
\max_{t\in[0,1]}\tilde I(\bar\gamma(t))< \tilde I(u_0).
\end{equation}
We sketch below the construction of such curve $\bar\gamma\in\Gamma$, see \cite[Lemma 4.3]{CN} for more details. 

\noindent $\bullet$ It is easy to see that there exist two positive numbers $t_-$ and $t_+$ ($t_-<1<t_+$), such that $t_-u_0\in U_-$ and $t_+ u_0\in U_+$.

\noindent $\bullet$ By ($f_3$), the function $t\in[t_-,t_+]\mapsto\tilde{I}(tu_0)$ has a unique strict maximum point at $t=1$. Hence, $$\tilde I(tu_0)<\tilde I(u_0)\quad\mbox{for all }t\in[t_-,t_+]\setminus\{1\}.$$
 
\noindent $\bullet$ Let $v\in W^{1,p}_\mathrm{rad}(B)\setminus\{0\}$ be nondecreasing and such that $\int_B vdx=0$. For every $t\in[t_-,t_+]$, the function $s\in\RR\mapsto\tilde I(t(u_0+sv))$ is continuous. Therefore, by the previous step, we get for $s$ in a neighborhood of $0$
$$\tilde I(t(u_0+sv))<\tilde I(u_0)\quad\mbox{for all }t\in[t_-,t_+]\setminus[1-\varepsilon,1+\varepsilon],$$
where $\varepsilon>0$ is a sufficiently small constant.

\noindent $\bullet$ In order to have the same inequality also for $t$ close to 1, we use condition ($f_3$), the $C^2$-regularity of $\tilde I$ and the Implicit Function Theorem, see \cite[Lemma 4.1]{CN}. This allows us to prove that $u_0$ is not a local minimum of the Nehari-type set 
$$\mathcal N_*:=\{u\in\mathcal C_*\,:\, \tilde I'(u)[u]=0\}.$$
In particular, we prove that for all $s\in\RR$ there exists a unique $\bar t_s>0$ such that $\bar t_s(u_0+sv)\in\mathcal N_*$ and $\bar t_s$ is the unique maximum point of the map  $t\in[1-\varepsilon,1+\varepsilon]\mapsto\tilde I(t(u_0+sv))$. Furthermore, by using a second order Taylor expansion of the energy functional and ($f_3$), we obtain that for $s$ in a neighborhood of $0$
$$\tilde I(\bar t_s(u_0+sv))-\tilde{I}(u_0)=\frac{s^2}{2}\int_B[(p-1)u_0^{p-2}-\tilde f'(u_0)]v^2dx+o(s^2)<0.$$
Therefore, we get for $s$ close to $0$
$$\tilde I(t(u_0+sv))\le \tilde I(\bar t_s(u_0+sv))<\tilde I(u_0)\quad\mbox{for all }t\in[1-\varepsilon,1+\varepsilon].$$ 

\noindent $\bullet$ Clearly, for $\bar s>0$ small enough, $t_-(u_0+\bar sv)\in U_-$ and $t_+(u_0+\bar s v)\in U_+$.

\noindent $\bullet$ By the convexity of $\mathcal C_*$, keeping in mind that $U_-,\, U_+\subset\mathcal C_*$,
$$t(u_0+\bar sv)\in\mathcal C_*\quad\mbox{for every }t\in[t_-,t_+].$$

\noindent $\bullet$ Hence, the curve $\bar \gamma:t\in[0,1]\mapsto ((1-t)t_-+tt_+)(u_0+\bar sv)\in\mathcal C_*$ belongs to $\Gamma$ and satisfies \eqref{toprove}.  
 
\section{Sketch of the proof of Theorem \ref{limit_profile}}
We denote by $u_q\in \mathcal C$ the nonconstant solution of 
\begin{equation}\label{Pq}
\begin{cases}
-\Delta_pu+u^{p-1}=u^{q-1}\quad&\mbox{in }B,\\
u>0&\mbox{in }B,\\
\partial_\nu u=0&\mbox{on }\partial B
\end{cases}
\end{equation}
at minimax level $c_q$ and by $\tilde I_q$ the energy functional associated to the corresponding truncated problem. We describe below the main steps to prove Theorem \ref{limit_profile}, see for reference \cite[Theorem 1.3]{CN} and also \cite{G}.

\noindent $\bullet$ In \cite[Lemma 5.5]{CN}, we find an a priori bound on $u_q$, uniform in $q$. Namely, 
$$\|u_q\|_{C^1(\overline{B})}\le C,\mbox{ with $C>0$ independent of $q\ge p+1$.}$$
Here we use the special form of $f$. 

\noindent$\bullet$ This ensures the existence of a limit profile $u_\infty$ for which
$$u_q\rightharpoonup u_\infty\mbox{ in }W^{1,p}(B)\quad\mbox{and}\quad u_q\to u_\infty\mbox{ in }C^{0,\mu}(\overline{B})\,\forall\,\mu\in(0,1)\quad\mbox{ as }q\to\infty.$$
Furthermore, $u_\infty(1)=1$, see \cite[Lemma 5.6]{CN}. 

\begin{remark}
By integrating over $B$ the first equation of problem \eqref{Pq}, we get $\int_B u_q^{p-1}(1-u_q^{q-p})dx=0$. Since $u_q>0$, $u_q\not\equiv 1$, and $u'_q\ge0$, it results 
$$u_q(0)<1\quad\mbox{and}\quad u_q(1)>1\quad\mbox{for all }\,q\ge p+1.$$
Heuristically, where $u_q\le \mathrm{Const.}<1$ (i.e., near the center of the ball $B$), $\lim_{q\to\infty}u_q^{q-1}=0$. So, it is natural to expect that $u_\infty$ solves $-\Delta_p u+u^{p-1}=0$ at least in a neighborhood of the origin. On the other hand, in the region where $u_q\ge1$ (i.e., in a neighborhood of $\partial B$), the same limit is an indeterminate form.
This is somehow responsible of the fact that the boundary condition is not preserved in the limit. We further remark that, by Hopf's lemma, $\partial_\nu G>0$ on $\partial B$, hence the $C^{0,\mu}(\overline{B})$-convergence is optimal.    
\end{remark}
\noindent$\bullet$ We introduce the quantity  
$c_\infty:=\inf\left\{\frac1p \|u\|_{W^{1,p}(B)}^p\,:\, u\in\mathcal C,\,u\big|_{\partial B}=1\right\}$
and we show that
$c_\infty=\inf\left\{\frac1p{\|u\|_{W^{1,p}(B)}^p}\,:\, u\in W^{1,p}(B),\,u\big|_{\partial B}=1\right\}$.
Furthermore, this infimum is {\it uniquely} achieved at $G$ (via the Direct Method of the Calculus of Variations), see \cite[Lemma 5.7]{CN}.

\noindent$\bullet$ We show in \cite[Lemma 5.8]{CN} that $c_\infty=\lim_{q\to\infty}c_q$. The proof relies mainly on the fact that the minimax level $c_q$ in the cone coincides with a Nehari-type level in the cone (also here we use the fact that $f$ is a pure power function), cf. \cite[Lemma 5.4]{CN}. As a consequence, we get that $c_\infty$ is attained at $u_\infty$ and $\|u_q\|_{W^{1,p}(B)}\to \|u_\infty\|_{W^{1,p}(B)}$. 

\noindent $\bullet$ By uniqueness, $u_\infty=G$ a.e. in $B$. Finally, the weak convergence ($u_q\rightharpoonup G$ in $W^{1,p}(B)$) together with the convergence of the norms ($\|u_q\|_{W^{1,p}(B)}\to \|G\|_{W^{1,p}(B)}$) guarantee that $u_q\to G$ in $W^{1,p}(B)$, by the uniform convexity of the space.
\medskip

\noindent{\bf Acknowledgements.}
The author gratefully thanks Dr. Benedetta Noris for her careful
reading of the manuscript and her valuable suggestions.

\begin{flushleft}

{\bf AMS Subject Classification: 35J92, 35B09, 35B40, 35B45, 35A15}\\[2ex]

Francesca~COLASUONNO,\\
D\'epartement de Math\'ematique, Universit\'e Libre de Bruxelles\\
Campus de la Plaine - CP214\\
boulevard du Triomphe, 1050 Bruxelles, Belgium\\
e-mail: \texttt{francesca.colasuonno@ulb.ac.be}\\[2ex]

\end{flushleft}

\end{document}